\documentclass[a4paper, 11pt]{article}

\usepackage{amsfonts,amsmath,amssymb,amscd,color,stmaryrd}
\usepackage[all]{xy} 
\usepackage{latexsym}
\usepackage[T1]{fontenc}

\newtheorem{thm}{Theorem}[section]
\newtheorem{definition}{Definition}[section]

\newtheorem{cor}[thm]{Corollary}

\newenvironment{prf}{\begin{trivlist}\item[]{\bf Proof: }}
{\hfill $\blacksquare$ \end{trivlist}}

\newcommand{\coker}{\text{coker}}
\newcommand{\ind}{\mathrm{index}\,}

\newcommand{\be}[1]{\begin{equation}\label{#1}}
\newcommand{\ee}{\end{equation}}

\newcommand{\beq}{\begin{eqnarray*}}
\newcommand{\eeq}{\end{eqnarray*}}
\newcommand{\bpr}{\begin{prf}}
\newcommand{\epr}{\end{prf}}

\newcommand{\Z}{\mathbb{Z}}

\newcommand{\mf}{\mathfrak}
\newcommand{\mb}{\mathbf}
\newcommand{\mc}{\mathcal}

\newcommand{\R}{\mathbb{R}}

\newcommand{\C}{\mathbb{C}}

\renewcommand{\O}{\mathbb{O}}
\renewcommand{\H}{\mathbb{H}}

\def\note#1{\marginpar{\raggedright\if@twoside\ifodd\c@page\raggedleft\fi\fi\sf\scriptsize RMK: #1}}

\def\ker{\mathop{\rm ker}\nolimits}

\def\Re{\mathop{\rm Re}\nolimits}

\def\Im{\mathop{\rm Im}}

\begin{document}
\title{Gauge theory in dimension $7$}
\author{Frederik Witt}
\date{}
\maketitle

\centerline{\textbf{Abstract}}
We first review the notion of a $G_2$--manifold, defined in terms of a principal $G_2$ (``gauge'') bundle over a $7$--dimensional manifold, before discussing their relation to supergravity. In a second thread, we focus on associative submanifolds and present their deformation theory. In particular, we elaborate on a deformation problem with coassociative boundary condition. Its space of infinitesimal deformations can be identified with the solution space of an elliptic equation whose index is given by a topological formula. 
 
\bigskip


\smallskip

\section{Introduction}
%%%%%%%%
%%%%%%%%
The rich geometric structures encountered in low dimensions stem from special algebraic features such as triality, which cease to be present in higher dimensions. Another thread is the existence of {\em vector cross products} $\times:\R^n\times\R^n\to\R^n$ which can only exist in dimension $n=3$ and $n=7$, and whose existence is tied to the algebraic structure of (imaginary) quaternions and (imaginary) octonions respectively. On a manifold $M$, the tangent bundle can be equipped with such a vector cross product if the structure (or ``gauge'') group $G$ reduces to the automorphism group of $\times$. If we also impose existence of a Riemannian metric $g$, we are left with $G=SO(3)$ or $G_2$~\cite{gr69}. This happens under fairly mild topological conditions: While for $n=3$ we only need orientability (i.e. the first Stiefel class of $M$ vanishes), in addition $n=7$ requires $M$ to be spin (i.e. the second Stiefel Whitney class vanishes). In this survey article, we focus on dimension $7$ and consider solely $G_2$--manifolds, though it is worthwhile to take dimension $3$ as a guidance~\cite{doth98}.

\bigskip

Since $G_2$ also appears on Berger's list, $G_2$--geometry has been investigated for a long time from the viewpoint of Riemannian holonomy, culminating with Joyce's celebrated construction of compact holonomy $G_2$--manifolds~\cite{jo96},~\cite{jo00}. Emphasis shifted when $G_2$--manifolds became important in supergravity compactifications. Here, the physical theory requires the holonomy to be contained in $G_2$ with respect to some connection $\widetilde{\nabla}$ which is not necessarily the Levi--Civita connection $\nabla^g$. In case $\widetilde{\nabla}$ is metric, the resulting condition can be regarded as a generalisation of Gray's concept of weak holonomy~\cite{gr71}. In particular, if $\widetilde{\nabla}$ is metric and torsion--free, $\widetilde{\nabla}$ coincides with $\nabla^g$, so that the underlying $G_2$--manifold is torsion--free and its holonomy is contained in $G_2$. We shall review these aspects as we go along in Section~\ref{g2geom}. 

\bigskip

A second line of thought is inspired by K\"ahler geometry. In real terms, K\"ahler manifolds are defined by a complex structure $J$ and a hermitian metric $g$, which give rise to the K\"ahler form $\omega(x,y)=g(Jx,y)$. On $G_2$--manifolds, the cross product $\times$ and the metric $g$ define the $3$--form $\varphi(x,y,z)=g(x\times y,z)$. Both forms are generic or {\em stable} following the language of Hitchin, and his variational principle puts torsion--free $G_2$-- and symplectic manifolds on equal footing~\cite{hi01}. Furthermore, $\omega$ and $\varphi$ distinguish special classes of submanifolds, namely complex and Lagrangian submanifolds in K\"ahler geometry and associative and coassociative submanifolds in $G_2$--geometry~\cite{hala82} (Section~\ref{strucman}). In the symplectic world, these submanifolds give rise to highly non--trivial invariants. Roughly speaking, the Gromov--Witten invariant counts the number of pseudo--holomorphic curves inside a symplectic manifold, while counting the number of pseudo--holomorphic Whitney discs bounding Lagrangians gives rise to Lagrangian intersection Floer homology. A first step towards the definition of similar invariants in the $G_2$--context is to study the deformation behaviour of associatives. For closed associatives, McLean showed that the deformation theory is governed by an elliptic equation whose index however is always $0$ on topological grounds~\cite{mc98}. Deformations of associatives with boundary inside a fixed coassociative were studied in~\cite{gawi08}. Again, the deformation problem gives rise to an elliptic equation. Its index is given by a topological formula for which examples with non--trivial index exist. These issues, as well as some technical aspects from PDE theory, will be discussed in Section~\ref{deftheo}.
%
%
%
%%%%%%%%
%%%%%%%%
\section{$G_2$--geometry}
%%%%%%%%
%%%%%%%%
\label{g2geom}
%
%%%%%%%%%%%
\subsection{The imaginary octonions}
%%%%%%%%%%%
\label{imoct}
In essence, $G_2$--geometry is the geometry of imaginary octonions. To fully appreciate this point of view, we shall discuss the complex counterpart of K\"ahler and Calabi--Yau geometry first. 

\bigskip

{\bf Hermitian spaces.}
In real terms, the structure of the standard hermitian space $(\C^m,h)$ is given by a complex structure $J$ on the underlying real vector space $V=\R^{2m}$, that is an endomorphism $J:V\to V$ squaring to minus the identity. We recover the complex space $\C^m$ as the $+i$--eigenspace $V^{1,0}$ of $J$ extended to the complexification $V\otimes\C=V^{1,0}\oplus\overline{V^{1,0}}$. Furthermore, $J$ is an isometry for the Euclidean inner product $g=\Re h$. We also say that $(J,g)$ defines a {\em K\"ahler structure} on $\R^{2m}$. Furthermore, we can define the {\em K\"ahler form} 
\begin{equation}\label{kaehlerform}
\omega(x,y)=g(Jx,y).
\end{equation}
Note that $GL(2m)$ acts both on the space of endomorphisms and positive definite Euclidean inner products in a natural way. The common stabiliser of $(J,g)$ is $U(m)$, and one therefore also refers to a K\"ahler structure as a $U(m)$--{\em structure}. 

\bigskip 

A special case of K\"ahler structures are Calabi--Yau structures which in terms of stabiliser groups are associated with $SU(m)\subset U(m)$. Apart from $(J,g)$ we are also given a {\em complex} volume form $\Omega\in\Lambda^mV^{1,0*}$ such that the two real $m$--forms $\psi_+=\Re\Omega$ and $\psi_-=\Im\Omega$ satisfy $\psi_+\wedge\psi_-=\omega^m/m!$.

\bigskip

{\bf Imaginary octonions.}
Next, consider the direct sum of two quaternionic spaces, namely the octonions $\O=\H\oplus e\H$ which is a real $8$--dimensional, non--associative division algebra generated by $\langle\mb{1},i,j,k,e,e\cdot i,e\cdot j,e\cdot k\rangle$. Taking these generators as an orthonormal basis induces an inner product $g$ on $\O$ compatible with the algebra structure. Further, we obtain a {\em cross product} $\times$ taking values in the {\em imaginary octonions} $\Im\O=\langle\mb{1}\rangle^{\perp}\cong\R^7$ by defining
$$
u\times v=\Im(\overline{v}\cdot u).
$$
Here, $\overline{v}$ is the natural conjugation which sends $v\in\Im\O$ to $-v$. The term cross product is justified by the properties $u\times v=-v\times u$ and $|u\times v|=|u\wedge v|$. In analogy to~(\ref{kaehlerform}), we can define the $3$--form
$$
\varphi_0(u,v,w)=g(u\times v,w),
$$
which expressed in the orthonormal basis $e_1=i,\,e_2=k,\ldots,e_7=e\cdot k$ can be written explicitly as
\begin{equation}\label{standform}
\varphi_0=e^{123}+e^1\wedge(e^{45}+e^{67})+e^2\wedge(e^{46}-e^{57})+e^3\wedge(-e^{47}-e^{56}).
\end{equation}
The stabiliser of $\varphi_0$ inside $GL(7)$ is $G_2$, which is why we refer to $\varphi_0$ as $G_2$--{\em form}. These exist in abundance: They are acted on transitively by $GL(7)$ so that the orbit of $G_2$--forms $GL(7)/G_2$ has dimension $49-14=35=\dim\Lambda^3\Im\O^*$ ($G_2$ being of dimension $14$). In particular, the orbit of $G_2$--forms is {\em open}. Further, $\varphi_0$ induces a volume form (which is somehow difficult to write down explicitly, cf. the appendix in~\cite{hi01}). This renders $G_2$--structures akin to Calabi--Yau structures, and in fact, starting from a Calabi--Yau structure on $\R^6$, the $3$--form
$$
\varphi_0=\psi_++\omega\wedge e^7
$$
induces a $G_2$--form on $\R^7=\R^6\oplus\R e_7$. On the level of stabiliser groups this is reflected by the inclusion $SU(3)\hookrightarrow G_2$, where $G_2/SU(3)=\mc{S}^6$ is the $6$--sphere in $\Im\O$.
%
%%%%%%%%%%%
\subsection{Topological and torsion--free $G_2$--manifolds}
%%%%%%%%%%%
\label{g2manif}
A $7$--dimensional manifold $M$ is called a {\em topological} $G_2$--{\em manifold} or simply a $G_2$--manifold if there exists $\varphi\in\Omega^3(M)$ such that $\varphi_x$ defines a $G_2$--structure on $T_xM$ as discussed in the previous section. By an {\em abus de langage}, we refer to the $3$--form $\varphi$ itself as the $G_2$--structure. This is tantamount to saying that the principal frame bundle associated with $GL(7)$ reduces to a $G_2$--principal frame bundle, which consists of isomorphisms between $(T_xM,\varphi_x)$ and $(\Im\O,\varphi_0)$ for $x\in M$. In particular, these isomorphisms induce a natural Riemannian metric $g$ in $M$. .

\bigskip

A $G_2$--structure is said to be {\em torsion--free} if $\nabla^g\varphi=0$, where $\nabla^g$ is the Levi--Civita connection associated with $g$. Equivalently there exist coordinates around each point such that $\varphi(x)=\varphi_0+O(|x|^2)$ so that the $G_2$--structure is flat to first order. The most important criterion for torsion--freeness is the theorem of Fern\'andez--Gray~\cite{fegr82}:

\begin{thm}\label{torsionfree}
A $G_2$--manifold $(M,\varphi)$ is torsion--free if and only if $d\varphi=0$ and $d\star\varphi=0$.
\end{thm}

\noindent The holonomy of a torsion--free $G_2$--metric is actually contained in $G_2$. In the sequel, we say that a torsion--free $G_2$--manifold is a {\em holonomy $G_2$--manifold}, if the holonomy equals $G_2$\footnote{Note that some authors do not make this distinction.}. 

\bigskip

A trivial example of a torsion--free $G_2$--structure is $\R^7$ with $\varphi$ as in~(\ref{standform}) (with the standard coordinates $dx^{ijk}$ in place of $e^{ijk}$). Since it is translation invariant, the $G_2$--structure descends to the torus $T^7=\R^7/\Z^7$ where it defines a compact torsion--free $G_2$--manifold. Examples of holonomy $G_2$--manifolds were constructed by Bryant~\cite{br87}, Bryant--Salamon, Joyce and Kovalev~\cite{ko03}. In~\cite{brsa89}, Bryant and Salamon define holonomy $G_2$--metrics on (an open set of) the total space of the spinor bundle $S\to M^3$, where $M^3$ is a three--dimensional space form. In particular, when $M$ is taken to be the $3$--sphere $\mc{S}^3$, there exists a {\em complete} holonomy $G_2$--metric on the total space $S\cong\mc{S}^3\times\H$ such that the fibres are orthogonal to the horizontal distribution of the canonical spin connection induced by $\nabla^g$. A method for the construction of compact holonomy $G_2$--manifolds is due to Joyce (\cite{jo96} and~\cite{jo00}). In essence, his idea consists in considering quotients $T^7/\Gamma$, where $\Gamma$ is a discrete group of isometries acting on $T^7$ which preserve the standard $G_2$--form $\varphi$. Therefore, $\varphi$ descends to a torsion--free $G_2$--form outside the singularity locus produced by dividing out the action of $\Gamma$. In favourable cases these can be resolved in such a way that the resolution $M\to T^7/\Gamma$ carries a $G_2$--structure $\widetilde{\varphi}$ with ``small'' torsion, which can then be deformed into a torsion--free $G_2$--structure by Joyce's deformation theorem. For instance, a suitable group $\Gamma$ is generated by
$$
\begin{array}{lcl}
\alpha(x_1,\ldots,x_7) & = & (x_1,x_2,x_3,-x_4,-x_5,-x_6,-x_7)\\
\beta(x_1,\ldots,x_7) & = & (x_1,-x_2,-x_3,x_4,x_5,\frac{1}{2}-x_6,-x_7)\\
\gamma(x_1,\ldots,x_7) & = & (-x_1,x_2,-x_3,x_4,\frac{1}{2}-x_5,x_6,\frac{1}{2}-x_7).
\end{array}
$$
The resulting $G_2$--structure has then holonomy $G_2$ on topological grounds. 
%
%%%%%%%%%%%
\subsection{$G_2$--manifolds in physics}
%%%%%%%%%%%
\label{g2physics}
To make contact with physics we have to give yet another characterisation of $G_2$--manifolds (cf. for instance~\cite{fkms97} or~\cite{wa89}). The physical literature on $G_2$--manifolds is extensive, and the list of references given below is by no means exhaustive.

\bigskip

{\bf Spinorial characterisation of $G_2$--manifolds.} 
As pointed out before, a $G_2$--manifold $(M,\varphi)$ carries a natural Riemannian metric and a volume form, or equivalently, an orientation. On a group level, this is tantamount to saying that $G_2\subset SO(7)$. Since $G_2$ is simply--connected, we can lift this inclusion to $Spin(7)$. Further, $Spin(7)/G_2\cong\mc{S}^7$, where $\mc{S}^7$ denotes the $7$--sphere in the real $8$--dimensional, irreducible spin representation of $Spin(7)$. Put differently, we can see $G_2$ not only as the stabiliser of a $3$--form of special algebraic type, but also as the stabiliser of a {\em unit spinor}. In global terms this means that the principal $G_2$--frame bundle induces a canonical spin structure with spinor bundle $S$. Further, $G_2$--manifolds carry a natural unit spinor field $\Psi\in C^{\infty}(M,S)$. Conversely, assume we are given a unit spinor field $\Psi$ for some spin structure on a $7$--dimensional Riemannian manifold $(M,g)$. Under the well--known identification $S\otimes S\cong\Lambda^*T^*M$, we have 
\begin{equation}\label{g2bispin}
\Psi\otimes\Psi=1+\varphi+\star\varphi+\mathrm{vol}_g 
\end{equation}
(cf. for instance~\cite{lami89} Section IV.10 or~\cite{wa89}). The difference between these two viewpoints is this: While the $G_2$--form is specified at each point by $35=\dim\Lambda^3\Im\O^*$ parameters, the spinor definition requires an a priori choice of a Riemannian metric $g$ which at each point is determined by $28=\dim\odot^2\Im\O^*$ parameters. The remaining $35-28=7=\dim\mc{S}^7$ degrees of freedom are fixed by the choice of a unit spinor field. By general principal fibre bundle theory, the $G_2$--structure defined in terms of $(g,\Psi)$ is torsion--free if and only if $\nabla^g\Psi=0$ holds, where by abus de notation, $\nabla^g$ denotes the Levi--Civita connection on the tangent bundle as well as the canonical lift to the spinor bundle.

\bigskip

{\bf Supersymmetry.} 
In physics, spinor field equations arise for instance in connection with supergravity and (super) string theory. Here is a rather informal explanation -- for the true and detailed story cf.~\cite{fr99}, or~\cite{ka06} for a shorter introduction. According to quantum mechanics there are two kinds of particles: {\em bosons} (which transmit forces such as photons) and {\em fermions} (which make up matter such as electrons). In the mathematical model building, bosons materialise as sections of tensor bundles (e.g. vector fields or differential forms) while fermions arise as sections of spinor bundles (e.g. spinor fields or spinor--valued differential forms). Now a {\em supersymmetry} is a symmetry taking fermions to bosons and vice versa, or, in more mathematical terms, a transformation from tensor bundles to spinor bundles. For instance, if we are given a spinor field $\Psi$, then Clifford multiplication induces a map taking vector fields $X\in C^{\infty}(M,TM)$ to spinor fields $X\cdot\Psi\in C^{\infty}(M,S)$. On physical grounds, one restricts attention to systems of fermions invariant under infinitesimal supersymmetry transformations (this is the so--called {\em localisation principle}) which leads to certain spinor field equations. We give two examples hereof next.

\bigskip

{\bf Heterotic supergravity and M--theory.} 
First we consider heterotic supergravity, the low energy limit of heterotic string theory, which takes place on a ten--dimensional Lorentzian spin manifold $N^{1,9}$. Supersymmetry materialises as before in terms of a unit spinor field $\Psi$. Furthermore, we have a $3$--form $H\in\Omega^3(N)$, the so--called {\em H--flux}. The localisation principle leads (among other constraints) to the {\em gravitino} equation 
\begin{equation}\label{hetsugra}
\nabla_X^g\Psi+\frac{1}{4}(X\llcorner H)\cdot\Psi=0.
\end{equation}
In order to solve this equation, one often makes a {\em compactification} ansatz of the form $\R^{1,p}\times M^p$ where $\R^{1,p}$ is now flat Minkowski space and $M^p$ a Riemannian manifold, usually taken to be compact (whence the name). In this case,~(\ref{hetsugra}) reduces to a spinor field equation on $M^p$ with $H\in\Omega^3(M)$. In particular, we obtain for $p=7$ a $G_2$--manifold $(M,g,\Psi)$. If we define the metric connection $\widetilde{\nabla}$ on $TM$ by
$$
g(\widetilde{\nabla}_XY,Z)=g(\nabla^g_XY,Z)+\frac{1}{2}H(X,Y,Z)
$$
for $X,\,Y,\,Z\in C^{\infty}(M,TM)$, then~(\ref{hetsugra}) is precisely the condition $\widetilde{\nabla}\Psi=0$, where we again abuse notation and denote by $\widetilde{\nabla}$ the natural lift to the spinor bundle. Geometrically speaking, this is just the assertion that the holonomy of $\widetilde{\nabla}$ is contained in $G_2$. If $H\equiv0$, then $\widetilde{\nabla}$ and $\nabla^g$ coincide, and we recover the condition for a torsion--free $G_2$--structure. In this sense, equation~(\ref{hetsugra}) can be seen as an extension of Gray's concept of {\em weak holonomy}~\cite{gr71}. A good mathematical reference is~\cite{friv01}, where Friedrich and Ivanov gave the first detailed account on this type of connections.

\bigskip

Another example is provided by M--theory. Here, we consider an eleven--dimensional Lorentzian spin manifold $N^{1,10}$ together with a unit spinor field and a $4$--{\em form flux} $F\in\Omega^4(N)$. Compactifying to $\R^{1,3}\times M^7$ as before yields a $G_2$--structure $(M,g,\Psi)$, where the spinor field $\Psi$ has to satisfy the equation
\begin{equation}\label{mtheo}
\widetilde{\nabla}_X\Psi=\nabla_X^g\Psi+\frac{1}{6}(X\llcorner F)\cdot\Psi+\frac{1}{12}(X\wedge F)\cdot\Psi=0
\end{equation}
(cf. for instance~\cite{beje03}). In contrast to the previous case, $\widetilde{\nabla}$ is not induced by a metric connection of $TM$, and understanding the geometric meaning for the underlying $G_2$--structure is less straightforward. Rather, one has to interpret this equation in terms of the holonomy of the spin bundle (leading to so--called {\em generalised holonomy} in physicists' jargon, cf. for instance~\cite{duli03}).
%
%
%
%%%%%%%%
%%%%%%%%
\section{Structured submanifolds}
%%%%%%%%
%%%%%%%%
\label{strucman}
In this section we introduce the notion of a calibrated submanifold as introduced by Harvey and Lawson in their seminal paper~\cite{hala82}. As they point out, an ambient geometric structure (say a complex manifold) can be investigated in terms of a distinguished family of submanifolds (say complex submanifolds). In the context of $G_2$--geometry, this eventually leads to the study of associative and coassociative submanifolds.
%
%%%%%%%%%%%
\subsection{Calibrations}
%%%%%%%%%%%
\label{calibrations}
{\bf Complex subspaces.} 
It is natural to ask whether there are any interesting substructures associated with $G_2$--geometry. Again, it is instructive to consider K\"ahler structures first. In $\C^m$ we have the natural notion of a complex subspace $V\subset\C^m$. In real terms, this means that the underlying real vector space $\llbracket U\rrbracket$ of $U$ is {\em stable} under the complex structure $J$, i.e. $J\big(\llbracket U\rrbracket\big)\subset\llbracket U\rrbracket$. 

\bigskip

{\bf Associative subspaces.} 
In $G_2$--geometry, the r\^ole of the complex structure $J$ is assumed by the cross product $\times$. A natural definition for a subspace $U\subset \Im\O$ is therefore to be stable under $\times$. The trivial dimensions $0$, $1$ and $7$ apart, a stable subspace is necessarily of dimension $3$. Harvey and Lawson call these subspaces {\em associative}, for stability is equivalent to the vanishing of the totally skew--symmetric {\em associator}
\begin{equation}\label{assoc}
[u,v,w]=\frac{1}{2}\big((u\cdot v)\cdot w-u\cdot(v\cdot w)\big).
\end{equation}
For example, the imaginary quaternions $\Im\H$ spanned by $i,\,j,\,k$ in the natural decomposition $\Im\O=\Im\H\oplus\H$ define an associative subspace. In fact, $G_2$ acts transitively on the set of associative subspaces, which is isomorphic to $G_2/SO(4)$~\cite{hala82}, so that associative spaces exist in abundance. Here, $SO(4)$ acts on $\H$ via its standard vector representation on $\R^4$, while the action on $\Im\H$ corresponds to one of the two non--trivial homomorphisms $\rho:SO(4)\to SO(3)$ (recall that $SO(4)\cong\big(Spin(3)\times Spin(3)\big)/\Z_2$). Then $A\in SO(4)$ acts on $\Im\H$ via $\rho(A)\oplus A$ as a subgroup of $G_2$. 

\bigskip

{\bf Calibrations.} 
More generally, Harvey and Lawson introduced calibrations to give a unified approach not only to complex and associative subspaces, but also to various natural substructures in further geometries. The general setting is given by a real (oriented) vector space  $(V,g,\tau)$ together with a Euclidean inner product $g$ and a $k$--form $\tau\in\Lambda^kV^*$. We say that $\tau$ defines a {\em calibration} if for every oriented $k$--subspace $\xi=e_1\wedge\ldots\wedge e_k$ in $V$ determined by some orthonormal oriented system $e_1,\ldots,e_k$, the inequality $\tau(e_1,\ldots,e_k)\le 1$ holds and is met for at least one $k$--plane. Such a plane is said to be {\em calibrated} by $\tau$. For example, the powers $\omega^m/m!$ of the K\"ahler form $\omega$ define a calibration, and the calibrated subspaces are precisely the complex subspaces with their natural orientation (of complex dimension $m$). In analogy to the K\"ahler case, (suitably oriented) associative subspaces are calibrated by $\varphi$, which is a direct consequence of the {\em associator equality}
$$
\varphi(x,y,z)^2+\frac{1}{4}|[x,y,z]|^2=|x\wedge y\wedge z|^2. 
$$
{\em Coassociative} subspaces are calibrated by the Hodge dual $\star\varphi$. Hence, they are perpendicular to associative subspaces and of dimension $4$. 
%
%%%%%%%%%%%
\subsection{Associative submanifolds}
%%%%%%%%%%%
{\bf Associatives and coassociatives.} 
Next let $(M,\varphi)$ be a $G_2$--manifold. The previous definition of (co--)associative subspaces gives a natural class of structured submanifolds for $M$:

\begin{definition}
A submanifold $Y$ of $M$ is said to be {\em associative} if $T_pY$ (regarded as a subspace of $\Im\O$ via a $G_2$--frame) is associative for all points $p\in Y$. Associative submanifolds are therefore necessarily of dimension $3$. Similarly, we say that a submanifold $X$ is {\em coassociative} if $T_pX$ is coassociative for all points $p\in X$. 
\end{definition}

\bigskip 

{\bf Calibrated submanifolds.} 
If an associative $Y$ is suitably oriented, it follows from the previous section that the $G_2$--form $\varphi$ restricts to the induced Riemannian volume form on $Y$, that is, associatives are {\em calibrated} in the sense of Harvey and Lawson. As a consequence of Stoke's theorem compact associatives are {\em absolute volume minimisers} in their homology class if the calibration form $\varphi$ is closed~\cite{hala82}. This is a far stronger condition than being {\em minimal} (vanishing mean curvature). Similarly, suitably oriented coassociatives are calibrated with respect to $\star\varphi$, and  homologically volume minimising if $\star\varphi$ is closed. 

\bigskip

{\bf Local equation.} 
While minimality of a submanifold is a second order condition, calibrations (inducing a first order condition) become a handy tool in finding minimal submanifolds. To construct examples, we first set out for finding associative submanifolds inside $\Im\O$. Since torsion--free $G_2$--manifolds are flat to first order (cf. Section~\ref{g2manif}), this will provide a quite reasonable local model for associative submanifolds, at least in the torsion--free case. Rather than testing the condition $\varphi(\xi)\equiv1$ for $\xi=x\wedge y\wedge z$, we test for the vanishing of the associator~(\ref{assoc}). We think of it as an $\Im\O$--valued $3$--form $\chi=(\chi^1,\ldots,\chi^7)^{\top}$, so that the condition on a $3$--submanifold $Y$ to be associative becomes $\chi_{|Y}\equiv0$ (as a matter of notation we denote here and in the sequel the pull--back of $\chi$ to $Y$ by $\chi_{|Y}$). The components $\chi^j$ generate algebraically a differential ideal $\mc{I}$ of $\Omega^*(\Im\O)$, whose $3$--dimensional integral manifolds inside $\Im\O$ are associative. Further, Cartan--K\"ahler theory can be invoked to show that every real analytic surface $\Sigma$ of $\Im\O$ (trivially integral as $\mc{I}$ is generated by forms of degree $3$) can be extended to a uniquely determined associative germ $Y$ containing $\Sigma$~\cite{hala82},~\cite{rosa07}. In fact, it follows from similar arguments that for every associative $E\subset T_pM$ of a torsion--free $G_2$--manifold $(M,\varphi)$, there exists an associative submanifold $Y\subset M$ with $T_pY=E$. 

\bigskip

The associativity condition has a beautiful reformulation as a partial differential equation involving the Dirac operator, which will serve as guidance for the deformation theory to be developed later. Let $f:U\subset\Im\H\to\H$ be a smooth function defined on some open domain $U$. Following~\cite{hala82}, the condition for $Y=\{\underline{x}\oplus f(\underline{x})\,|\,\underline{x}\in\Im\H\}\subset\Im\H\oplus\H=\Im\O$ to be associative is this:

\begin{thm}
Let $f:U\subset\Im\H\to\H$ be a smooth function. Then $Y=\mathrm{graph}f$ is associative if and only if
$$
D(f)=i\frac{\partial f}{\partial x_1}+j\frac{\partial f}{\partial x_2}-k\frac{\partial f}{\partial x_3}=\sigma(\frac{\partial f}{\partial x_1},\frac{\partial f}{\partial x_2},\frac{\partial f}{\partial x_3}),
$$
where $D$ is the Dirac operator\footnote{The minus sign in front of the $k$ is due to our conventions which are based on~\cite{jo04}.} on $\Im\H$ and $\sigma:\H\times\H\times\H\to\H$ is the so--called triple cross product on $\H$.
\end{thm}

\bigskip

{\bf Closed examples.} In general, any closed, real analytic Riemannian $3$--manifold can be isometrically embedded as an associative into some (in general incomplete) torsion--free $G_2$--manifold~\cite{rosa07}. Further, consider the {\em complete} Bryant--Salamon metric on the total space of the spinor bundle $S\to\mc{S}^3$ over the $3$--sphere $\mc{S}^3$ (cf. Section~\ref{g2manif}). Here, the zero section $\mc{S}^3\times\{0\}$ defines an associative. Trivial compact examples are provided by extending compact Calabi--Yau $3$--folds $(K,\omega,\Omega)$ to $G_2$--manifolds via $(M=K\times S^1,\varphi=\Re\Omega+\omega\wedge dt)$ (cf. also Section~\ref{imoct}). A {\em complex curve} $C\subset K$ induces then the associative $C\times \mc{S}^1$. 

\bigskip

A general method for finding associatives inside a torsion--free $G_2$--manifold $(M,\varphi)$ is due to Joyce (cf. 10.8 in~\cite{jo00}). Let $\sigma:M\to M$ be an isometric, non--trivial involution on a torsion--free $G_2$--manifold $(M,\varphi)$ such that $\sigma^*\varphi=\varphi$. Then the fixed point locus of $\sigma$ defines an associative. This becomes a practical tool for the construction of associatives inside Joyce manifolds. For instance, consider the isometric involution $\sigma_0:T^7\to T^7$ given by
$$
\sigma_0(x_1,\ldots,x_7) = (x_1,x_2,x_3,\frac{1}{2}-x_4,-x_5,-x_6,-x_7).
$$
It satisfies $\sigma_0^*\varphi_0=\varphi_0$ and descends to an isometric involution $\sigma$ on $T^7/\Gamma$ with $\Gamma$ as given in Section~\ref{g2manif}, for $\sigma_0$ commutes with $\Gamma$. Then one can resolve $T^7/\Gamma$ in a $\sigma$--{\em equivariant} way, that is $\sigma$ lifts to an isometric involution on the resolution $(M,\varphi)$ such that $\sigma^*\varphi=\varphi$. The fixed point locus of $\sigma$ therefore defines an associative inside $M$.

\bigskip

{\bf Calibrations in physics.} In string and M--theory, {\em branes} are extended objects which minimise a certain energy functional. In the most simple cases, branes can be thought of as minimal submanifolds. But there is more to it -- namely constraints coming from supersymmetry which tell us that branes are not merely minimal, but calibrated. For this one needs to relate spinors with calibrations, which has been worked out by Dadok and Harvey~\cite{daha93},~\cite{ha90}. For instance, consider a $G_2$--manifold $(M,g,\Psi)$. Then the homogeneous components of the bi--spinor~(\ref{g2bispin}) define calibration forms. Further, a submanifold $Y$ with Riemannian volume form $\mathrm{vol}_Y$ acts on spinor fields via Clifford multiplication, and $Y$ is calibrated precisely if $\mathrm{vol}_Y\cdot\Psi=\Psi$ holds. On the other hand, if $\Psi$ is a spinor field parallel with respect to the modified spin connection $\widetilde{\nabla}$ (cf. Section~\ref{g2physics}), then this is the condition on $Y$ to represent a supersymmetric brane (cf. for instance~\cite{ga03} Section 4).

\bigskip

As discussed in Section~\ref{g2physics}, the metric spin connection $\widetilde{\nabla}$ one considers in supergravity will usually have torsion, which prevents the calibration forms from being closed as can be seen from Theorem~\ref{torsionfree}. Therefore, the calibrated submanifolds are not necessarily volume minimising. Rather, they minimise the {\em (brane) energy} $\mc{E}(X)=\mathrm{vol}(X)-\int_X\gamma$, where $d\tau=d\gamma$~\cite{gip03},~\cite{gupa99}. The form $\gamma$ can be interpreted as {\em Ramond--Ramond potential}, and $\int_X\gamma$ as the {\em Wess--Zumino term} of the brane energy. For a further development of these ideas in the context of so--called {\em generalised geometries}, see~\cite{gmwi08}.
%
%
%
%%%%%%%%
%%%%%%%%
\section{Deformations}
%%%%%%%%
%%%%%%%%
\label{deftheo}
Let $Y$ be a structured submanifold of some ambient geometry, for instance a complex submanifold inside a K\"ahler manifold or an associative submanifold inside a $G_2$--manifold. A natural object of study is the {\em moduli space} $\mf{M}_Y$ of all structured submanifolds isotopic to $Y$. A basic problem is to determine the Zariski tangent space of $\mf{M}_Y$, that is, the space of first order deformations of $Y$.
%
%%%%%%%%%%%
\subsection{Closed associatives}
%%%%%%%%%%%
{\bf Deformation of closed coassociatives.} 
Though we are primarily interested in the deformation theory of associatives, for motivating the later development it is instructive to start with the coassociative case first. The central result is due to McLean~\cite{mc98}:

\begin{thm}
Let $X$ be a closed coassociative (i.e. compact and without boundary) inside a $G_2$--manifold $(M,\varphi)$ with $d\varphi=0$. Then $\mf{M}_X$ is a smooth manifold of dimension $b^2_+(X)$, the dimension of real positive--definite $2$--cohomology.
\end{thm}

\noindent Let us briefly sketch the techniques of the proof which are quite archetypical (see also~\cite{jo04}). First we try to describe the set of nearby coassociatives by a smooth equation. To that end we fix a tubular neighbourhood $\mc{U}$ of $X$ which we think of as an open subset of the normal bundle $\nu\to X$ around the zero section. Submanifolds $X'$ which are $C^1$--close to $X$ correspond then to sections of $\mc{U}$ under exponentiation. For $X'$ to be coassociative we need $\varphi_{|X'}\equiv0$. Since $X'$ is isotopic to $X$ we can pull back $\varphi_{|X'}$ to $X$, where it lies in the same cohomology class as $[\varphi_{|X}]=0$, that is, the pull--back is {\em exact}. We obtain thus a smooth map $F:C^{\infty}(X,\mc{U})\to B^3(X)$ for which $F^{-1}(0)$ consists precisely of the coassociatives close to $X$. 

\medskip

Next we determine the space of first order deformations of $X$, that is, the kernel of the linearisation $d_0F$ of $F$ at the zero section. This is the so--called {\em Zariski tangent space} $T^{\mathrm{Zar}}\mf{M}_X\subset C^{\infty}(X,\nu)$ of $\mf{M}_X$. It consists of normal vector fields $s$ with $(\mc{L}_s\varphi)_{|X}=0$ , where $\mc{L}_s$ denotes the Lie derivative along $s$. In fact, one can show that $s\in C^{\infty}(X,\nu)\mapsto (s\llcorner\varphi)_{|X}\in\Omega^2(X)$ induces a bundle isomorphism between $\nu$ and the bundle of self--dual $2$--forms $\Lambda^2_+X$ of $X$. Since $(\mc{L}_s\varphi)_{|X}=d(s\llcorner\varphi)_{|X}$ (here we use the assumption $d\varphi=0$), the Zariski tangent space becomes the space of closed (and therefore coclosed, i.e. harmonic) $2$--forms under this identification. By standard Hodge theory, the dimension of this space is $b_+^2(X)=\dim\mc{H}^2_+(X)$. In particular, $T^{\mathrm{Zar}}\mf{M}_X$ can be regarded as the solution space of an {\em elliptic} equation. 

\medskip

In general, arbitrary first order deformations will not be realised as the deformation vector field of an actual deformation, which is why the dimension of the Zariski tangent space is sometimes referred to as the {\em virtual dimension} of the moduli space. In the present case however, $\mf{M}_X$ is smooth. For this, we call on the following version of the {\em Implicit Function Theorem} (see~\cite{amr88} Section 2.5 for this and related variations on that theme):
\begin{quote}
{\em Let $F:U\subset V\to W$ be a smooth map from some open neighbourhood $U$ around the origin of a Banach space $V$ into some other Banach space $W$.  If $\ker d_0F$ is finite--dimensional and $d_0F:V\to W$ is surjective, then the fibre $F^{-1}(0)$ is a smooth manifold locally isomorphic to $\ker d_0F$.}
\end{quote}
Now the spaces $C^{\infty}(X,\nu)$ and $B^3(X)$ are not Banach, so we extend $F$ to a smooth map $F^{k,\gamma}$ from $C^{k+1,\gamma}(X,\mc{U})$ inside a suitable H\"older space $C^{k+1,\gamma}(X,\Lambda^2_+X)$, $k\geq1$, $\gamma\in(0,1)$, to the Banach subspace of exact $C^{k,\gamma}$ $3$--forms inside $C^{k,\gamma}(X,\Lambda^3T^*X)$. One can then verify the surjectivity of $d_0F^{k,\gamma}$. Further, $F^{k,\gamma}=0$ is still an elliptic equation, so not only is the kernel finite--dimensional, but consists of {\em smooth} sections, that is, $\ker d_0F^{k,\gamma}=\ker d_0F$. Hence, $\mf{M}_X$ is a smooth manifold locally isomorphic to the space of harmonic self--dual $2$--forms.

\medskip

As an example, take $X$ to be a K3 surface $K$ or a $4$--torus $T^4$.  Both are real analytic Riemannian manifolds whose bundle of self--dual $2$--forms $\Lambda^2_+X$ is trivial. By a theorem of Bryant's~\cite{br00}, they can be isometrically embedded into a torsion--free $G_2$--manifold. Since in both cases $b^2_+(X)=3$, $X$ moves in a $3$--dimensional coassociative family. Actually, $X$ can be embedded as the $0$--fibre of a fibration $M^7\to B^3$ with coassociative fibres, where $B^3$ is a neighbourhood of $0\in\R^3$. This is reminiscent of the SYZ--formulation of Mirror Symmetry~\cite{syz96} which involves Calabi--Yau $3$--folds fibred by special Lagrangians, and indeed, there are corresponding conjectures for coassociative fibrations of torsion--free $G_2$--manifolds in connection with M--theory~\cite{gyz02}.

\bigskip

{\bf Deformations of associatives.} 
Next we address deformations of a closed associative $Y$ inside a $G_2$--manifold $(M,\varphi)$. Let us start with some heuristic considerations. The set of associative subspaces in $\Im\O$ is diffeomorphic to $G_2/SO(4)$, a codimension $4$ submanifold in the Grassmannian $G_3(\Im\O)$ of $3$--planes in $\Im\O$ (cf. Section~\ref{calibrations}). The condition on a $3$--plane to be associative is therefore (locally) given by four independent equations. On the other hand, we are free to vary along the normal bundle $\nu\to Y$ which now is of rank $4$, so the deformation problem involves four equations in four functions. It is therefore a {\em determined} problem. To get a feeling for the Zariski tangent space, we consider deformations of $Y=\Im\H\oplus0\subset\Im\O$ which we think of as the graph of $f_0\equiv0$. Nearby deformations are given by the graphs of a smooth family of functions $f_t:\Im\H\to\H$, and the deformation vector field is the partial derivative $s(\underline{x})=\frac{\partial}{\partial t}f_t(\underline{x})_{|t=0}$. As we have seen before, $Y_t=\mathrm{graph}\,f_t$ is associative if and only if $\mb{D}(f)=(D-\sigma)(f)=0$. Now $\epsilon^{-1}\mb{D}(\epsilon s)=D(s)-\epsilon^2\sigma(s)$, so that the linearisation $d_0\mb{D}(s)=\lim_{\epsilon\to0}\mb{D}(\epsilon s)/\epsilon=D(s)$ characterises the Zariski tangent space as the solution space of the {\em Dirac equation} $D(s)=0$. 

\medskip

The previous considerations generalise as follows. Let $Y$ be an associative inside some $G_2$--manifold $(M,\varphi)$. One can identify the normal bundle $\nu\to Y$ with a {\em twisted} spinor bundle, and choosing a connection $\nabla$ on $\nu$ induces an associated Dirac operator $\mb{D}^{\nabla}:C^{\infty}(Y,\nu)\to C^{\infty}(Y,\nu)$. The Zariski tangent space is then characterised by the following generalisation of McLean's theorem due to Gutowski, Ivanov and Papadopoulos~\cite{gip03} and Akbulut and Salur~\cite{aksa04},~\cite{aksa07}:

\begin{thm}\label{mclean}
Let $Y$ be a closed associative inside a $G_2$--manifold $(M,\varphi)$. Then there exists a connection $\nabla$ on $\nu$ such that
$$
T^{\mathrm{Zar}}\mf{M}_Y\cong\ker\mb{D}^{\nabla},
$$
where $\mb{D}^{\nabla}$ denotes the Dirac operator associated with $\nabla$.
\end{thm}

\noindent If the $G_2$--manifold is torsion--free, the theorem holds for the natural connection on $\nu$ induced by the Levi--Civita connection of $M$, and we recover McLean's original result as proven in~\cite{mc98}.

\medskip

Again, the Zariski tangent space is the solution space to an elliptic equation for which the index $\ind(\mb{D}^{\nabla})=\dim\ker\mb{D}^{\nabla}-\dim\coker\mb{D}^{\nabla}$ is defined. Since $Y$ is odd--dimensional, $\ind(\mb{D}^{\nabla})=0$. In a generic situation, where one expects the cokernel to vanish, the virtual dimension would be zero as a consequence. In this sense, associatives are {\em virtually rigid}, and by counting these, one could hope to define an invariant of the underlying $G_2$--structure in analogy to Gromov--Witten invariants.
%
%%%%%%%%%%%
\subsection{Associatives with boundary}
%%%%%%%%%%%
We are now going to consider deformation problems with boundary (see also~\cite{kolo07} for a boundary problem in some sense inverse to ours). As a result, we will be able to derive a topological formula for the virtual dimension of the moduli space we consider.

\bigskip

{\bf Coassociative boundary condition.} Let $(M,\varphi)$ be a $G_2$--manifold, $X\subset M$ a coassociative and $Y\subset M$ a compact associative with boundary $\partial Y\subset X$. We wish to investigate the moduli space
$$
\mf{M}_{X,Y}=\{Y'\,|\,Y'\mbox{ compact associative isotopic to }Y\mbox{ with }\partial Y'\subset X\}.
$$
As for the closed case, we need to analyse the normal bundle $\nu\to Y$ first. Apart from being a twisted spinor bundle, more can be said near the boundary. To that end fix a collar neighbourhood $\mc{C}\cong\partial Y\times[0,\epsilon)$ of $\partial Y$ inside $Y$. Let $u$ denote the inward pointing unit vector field defined on $\mc{C}$. It follows from the properties of the cross product $\times$ (cf. Section~\ref{imoct}) that $u$ induces a hermitian structure near the boundary, namely
$$
G:\nu\to\nu,\quad G(x)=u\times x.
$$
This acts indeed as an isometry with respect to $g$, as
$$
g(Ga,Gb)=\varphi(u,a,u\times b)=-g\big(u\times(u\times b),a\big)=g(a,b)
$$
for any $a,\,b\in\nu_{|\mc{C}}$. Let $\nu_X\subset TX_{|\partial Y}$ denote the orthogonal complement of $T\partial Y $ in $TX_{|\partial Y}$. Then~\cite{gawi08} 
\begin{itemize}
	\item the bundle $\nu_{X}$ is contained in $\nu$ and is stable under $G$,
	\item the orthogonal complement $\mu_X$ of $\nu_X$ in $\nu$ is also stable under $G$, and
	\item viewing  $T\partial Y$, $\nu_X$ and $\mu_X$ as $G$--complex bundles, we have 
\begin{equation}\label{muvsnu}
\overline{\mu}_X\cong\nu_X\otimes_{\C}T\partial Y,
\end{equation}
that is $\mu_X^{0,1}\cong\nu_X^{1,0}\otimes T^{1,0}\partial Y\cong\nu^{1,0}_X\otimes\overline{K}_{\partial Y}$, where $K_{\partial Y}$ is the canonical line bundle over $\partial Y$.
\end{itemize}
Consequently, as the deformation vector field $s$ of a curve $Y_t\subset\mf{M}_{X,Y}$ has to be tangent to $X$, we must have $s_{|\partial Y}\in C^{\infty}(\partial Y,\nu_X)$. So, if we let 
\begin{equation}\label{Bop}
\mb{B}:C^{\infty}(Y,\nu)\to C^{\infty}(\partial Y,\mu_X)
\end{equation}
be the real operator of order $0$ which projects smooth sections of $\nu$ to $\mu_X$ over $\partial Y$, then as a corollary to (the generalised version of) McLean's theorem (Theorem~\ref{mclean}), the Zariski tangent space of $\mf{M}_{X,Y}$ is given by
$$
T_{\mathrm{Zar}}\mf{M}_{X,Y}\cong\ker\left(\begin{array}{lcc}
\mb{D}& & C^{\infty}(Y,\nu)\\
\oplus:C^{\infty}(Y,\nu)&\to&\oplus\\
\mb{B}& &C^{\infty}(\partial Y,\mu_X)
\end{array}\right).
$$

\bigskip

{\bf Boundary problems for Dirac operators.} Again, we would like to compute the virtual dimension as the index of the differential operator $\mb{D}\oplus\mb{B}$. This requires a suitable notion of ellipticity for this problem. In particular, we demand the following two properties:
\begin{itemize}
	\item {\em Fredholm property} for $\mb{D}\oplus\mb{B}$: the kernel and cokernel are finite dimensional, so the {\em index} $\ind(\mb{D}\oplus\mb{B})=\dim(\ker\mb{D}\oplus\mb{B})-\dim(\coker\,\mb{D}\oplus\mb{B})$ is defined.
	\item {\em Regularity}: if $f\in H^s(Y,\nu)\cap\ker(\mb{D}\oplus\mb{B})$, then $f\in C^{\infty}(Y,\nu)$ (where $H^s(Y,\nu)$, $s\geq0$ denotes the standard chain of Sobolev spaces; in particular, we have $H^0(Y,\nu)=L^2(Y,\nu)$, the square integrable sections of $\nu$).
\end{itemize}
Before we can define a suitable elliptic boundary condition, we need to introduce the {\em Calder\'on projector} $\mc{Q}_{\mb{D}}$ associated with the Dirac operator\footnote{As emphasised in the introduction of~\cite{bowo93}, the subsequent statements hold for any operator $\mb{D}$ of {\em Dirac type}, that is, the principal symbol of $\mb{D}^2$ satisfies $\sigma(\mb{D}^2)(x,\xi)=||\xi||^2$.} $\mb{D}$ (cf.~\cite{bowo93}, Thm. 12.4). This is a pseudo--differential operator 
$$
\mc{Q}_{\mb{D}}:C^{\infty}(\partial Y,\nu)\to\mc{C}(\mb{D})=\{s_{|\partial Y}\in C^{\infty}(\partial Y,\nu)\,|\,s\in C^{\infty}(Y,\nu),\,\mb{D}s=0\}
$$ 
of order $0$ mapping the smooth sections of $\nu$ over $\partial Y$ to the space of Cauchy data of\footnote{We are glossing over some technical details such as the passing to the ``closed double'' $M=Y\cup_{\partial Y}Y$, cf. Chapters 9, 11 and 12 in~\cite{bowo93}.} $\mb{D}$. Let $q=\sigma(\mc{Q}_{\mb{D}})$ denote the principal symbol of $\mc{Q}_{\mb{D}}$, which becomes important in the following

\begin{definition} {\rm (cf.~\cite{bowo93} Def. 18.1)} 
Let $Y$ be an arbitrary smooth manifold with boundary and $\nu\to Y$ a (twisted) spinor bundle. A  pseudo--differential operator $\mb{B}:C^{\infty}(\partial Y,\nu)\to C^{\infty}(\partial Y,\nu)$ of order $0$ is said to define an {\em elliptic boundary condition} (abbreviated e.b.c.) if and only if
\begin{itemize}
	\item the extension $\mb{B}^{(s)}:H^s(\partial Y,\nu)\to H^s(\partial Y,\nu)$ has closed range.
	\item the restriction of the principal symbol $b=\sigma(\mb{B})_{|\mathrm{range}(q)}:\mathrm{range}(q)\to\mathrm{range}(b)$ is an isomorphism. 
\end{itemize}
\end{definition}

\noindent If $\mb{B}$ defines an e.b.c., then regularity holds (\cite{bowo93} Thm. 19.1).  An example of an e.b.c. is the Atiyah--Patodi--Singer boundary condition~\cite{aps75}. In our situation, a more stringent condition holds:

\begin{definition} {\rm (cf.~\cite{bowo93} Rem. 18.2)}
An e.b.c. is said to be {\em local}, if in addition $\mathrm{range}(p,\xi)=\nu_p$ holds for all $p\in\partial Y$. 
\end{definition}

\noindent For a local e.b.c. $\mb{D}\oplus\mb{B}$ is a Fredholm operator whose index is given by the index of a Fredholm operator on the boundary, namely
$$
\ind(\mb{D}\oplus\mb{B})=\ind\big(\mb{B}\mc{Q}_{\mb{D}}:\mc{C}(\mb{D})\to C^{\infty}(\partial Y,\nu)\big)
$$
(\cite{bowo93} Thm. 20.12). Furthermore, this integer depends only on the homotopy type of the principal symbols involved (~\cite{bowo93} Thm. 20.13 and Rem. 22.25). We note that for even--dimensional manifolds the existence local e.b.c. is topologically obstructed (\cite{bobl85} Section II.7.B). For odd--dimensional manifolds (as in the case of an associative $Y$), the orthogonal projector\footnote{By an {\em orthogonal projector} we understand an operator $\mb{P}$ of order $0$ satisfying $\mb{P}=\mb{P}^2=\mb{P}^*$.} $\mb{P}^+$ onto $\nu^+$, the bundle of positive half--spinors over $\partial Y$, defines a local e.b.c. with vanishing index. Furthermore, the difference between the index of two local e.b.c. $\mb{B}_{1,2}:C^{\infty}(Y,\nu)\to C^{\infty}(\partial Y,\nu_{1,2})$ (where $\nu_{1,2}$ can be bundles different from $\nu$) is the index of a Fredholm operator over the boundary, namely
\begin{equation}~\label{diffind}
\ind(\mb{D}\oplus\mb{B}_2)-\ind(\mb{D}\oplus\mb{B}_1)=\ind\big(\mb{B}_2\mc{Q}_{\mb{D}}\mb{B}_1^*:C^{\infty}(\partial Y,\nu_1)\to C^{\infty}(\partial Y,\nu_2)\big)
\end{equation}
(\cite{bowo93} Thm. 21.2). With these tools at hand, one is in a position to prove~\cite{gawi08}

\begin{thm}\label{defthm}
The operator $\mb{B}:C^{\infty}(\partial Y,\nu)\to C^{\infty}(\partial Y,\mu_X)$ as defined in~(\ref{Bop}) induces a local e.b.c.. Furthermore,
$$
\ind(\mb{D}\oplus\mb{B})=\ind(\overline{\partial}_{\nu_X}),
$$
where $\overline{\partial}_{\nu_X}$ denotes the Cauchy--Riemann operator of $\nu_X$ (regarded as a complex line bundle).
\end{thm}

\begin{cor} 
If the boundary is connected, the Riemann--Roch theorem yields 
$$
\ind(\mb{D}\oplus\mb{B})=\int_{\partial Y}c_1(\nu_X)+1-g,
$$ 
where $g$ is the genus of $\partial Y$ and $c_1(\nu_X)$ is the first Chern class of $\nu_X$ with respect to the natural complex structure induced by $u$. 
\end{cor}

\noindent For the proof of Theorem~\ref{defthm} we first note that in a collar neighbourhood of the boundary, we can write $\mb{D}=u\cdot(\partial_u+\mb{C})$, where $\mb{C}$ is the so--called {\em tangential part} of $\mb{D}$. The principal symbol of $\mc{Q}_{\mb{D}}$ can then be computed from the principal symbol of $\mb{C}$ (\cite{bowo93} Thm. 12.4), and one can check that the condition for a local e.b.c. holds. Furthermore,~(\ref{diffind}) implies (with $\mb{B}_2=\mb{P}^+$ whose index is zero as remarked before) the index of $\mb{D}\oplus\mb{B}$ to be the index of an operator over the boundary. Making use of~(\ref{muvsnu}), one finally shows that its principal symbol coincides with the principal symbol of the Cauchy--Riemann operator $\overline{\partial}_{\nu_X}$.

\bigskip

{\bf An example.} 
We conclude with an example of non--zero index. 
%Further examples, including associatives in compact holonomy $G_2$--manifolds constructed via Joyce's method, can be found in~\cite{gawi08}. 
Let $(M,\varphi)$ be a torsion--free $G_2$--manifold, and $Y$ an associative submanifold $Y$  with real analytic boundary $\partial Y$. For instance, take $M=\R^7$ and $\partial Y$ a compact oriented Riemann surface of genus $g$. Let $a\in C^{\infty}(\partial Y,\nu)$ be a nowhere vanishing real analytic section. Since the metric of a torsion--free $G_2$--manifold is necessarily Ricci flat~\cite{bo66}, the metric is real analytic in harmonic coordinates~\cite{dk81}. Consequently, so is the geodesic flow $\gamma_a:\partial Y\times(-\epsilon,\epsilon)\to M$ induced by $a$, which therefore generates an analytical submanifold $N$ of dimension $3$. Further, $\varphi(v,w,a)=0$ for $v,w\in T\partial Y$, and since $\nabla\varphi=0$, we conclude that the pull--back of $\varphi$ to $N$ vanishes identically. A Cartan--K\"ahler type argument invoked by Harvey and Lawson~\cite{hala82} (see also~\cite{br00}) shows that $N$ determines a unique coassociative germ $X$ containing $N$. Furthermore, $\nu_X$ is generated by $a$ and $u\times a$, where $u$ denotes again the inward pointing normal vector field of $\partial Y$. Hence $c_1(\nu_X)=0$ and therefore the index equals $1-g$.

\bigskip

\begin{center}
{\bf Acknowledgments}
\end{center}
I wish to warmly thank Damien Gayet for the fruitful collaboration and the organisers of the very enjoyable workshop ``Mathematicial Challenges in String Phenomenology'' at ESI, Vienna, for the invitation and the hospitality provided whilst writing this review.

\bigskip

\noindent
F.~\textsc{Witt}: NWF I -- Mathematik, Universit\"at Regensburg, D--93040 Regensburg, F.R.G.\\
e-mail: \texttt{frederik.witt@mathematik.uni-regensburg.de}
\end{document}